\documentclass[12]{amsart}
\usepackage{amssymb,amsfonts,latexsym}

\newtheorem{theorem}{Theorem}[section]
\newtheorem{lemma}[theorem]{Lemma}

\newtheorem{corollary}[theorem]{Corollary}

\theoremstyle{definition}
\newtheorem{definition}[theorem]{Definition}

\newtheorem{conjecture}[theorem]{Conjecture}
\newtheorem{remark}[theorem]{Remark}
\newtheorem{general remarks}[theorem]{General remarks}

\newcommand{\id}{\operatorname{id}}

\renewcommand{\span}{\operatorname{span}}

\newcommand{\ben}{\begin{enumerate}}
\newcommand{\een}{\end{enumerate}}

\begin{document}
\title[On the codimension growth of G-graded algebras] {On the
codimension growth of $G$-graded algebras}

\author{Eli Aljadeff}
\address{Department of Mathematics, Technion-Israel Institute of
Technology, Haifa 32000, Israel}
\email{aljadeff@tx.technion.ac.il}

\date{Nov. 10, 2009}



\keywords{graded algebra, polynomial identity}

\thanks {The author was partially supported by the ISRAEL SCIENCE FOUNDATION
(grant No. 1283/08) and by the E.SCHAVER RESEARCH FUND}


\begin{abstract} Let $W$ be an associative \textit{PI} affine
algebra over a field $F$ of characteristic zero. Suppose $W$ is
$G$-graded where $G$ is a finite group. Let $\exp(W)$ and
$\exp(W_{e})$ denote the codimension growth of $W$ and of the
identity component $W_{e}$, respectively. We prove: $\exp(W)\leq
|G|^2 \exp(W_{e}).$ This inequality had been conjectured by Bahturin
and Zaicev.
\end{abstract}

\maketitle

\begin{section}{Introduction} \label{Introduction}

Let $W$ be a \textit{PI}-affine algebra over a field $F$ of
characteristic zero. The codimension growth of $W$ was studied by
several authors (see for instance \cite{BZ},
\cite{benanti-giamb-pipitone}, \cite{GMZ}, \cite{GMZ2}, \cite{gz1},
\cite{gz2}, \cite{GZbook}, \cite{reg2}, \cite{Z}). It provides an
important tool which measures the ``size" of the $T$-ideal of
identities of $W$ in asymptotic terms. In case the algebra~$W$ is
$G$-graded, it is natural to compare the codimension growths of $W$
and $W_{e}$, where $W_{e}$ is the identity component of $W$ with
respect to the given $G$-grading.

Let us recall briefly the definitions. Let $X$ be a countable set of
indeterminates and $F\langle X\rangle$ the corresponding free
algebra over $F$. Let $\id(W)$ be the $T$-ideal of identities of $W$
in $F\langle X\rangle$ and let $\mathcal{W}= F\langle
X\rangle/\id(W)$ denote the relatively free algebra of $W$. We denote by $c_n(W)$
the dimension of the subspace spanned by
multilinear elements in $n$ free generators in $\mathcal{W}$. Our
interest is in the asymptotic behavior of the sequence of
codimensions, namely in
$$
\exp(W)=\lim_{n\to \infty}\sqrt[n]{c_n(W)}.
$$

It is known that $\exp(W)$ exists and moreover it assumes only
integer values (see~\cite{gz1}, ~\cite{gz2}).

Now suppose that the algebra~$W$ is $G$-graded, where $G$ is an
arbitrary finite group. Clearly, the identity component $W_e$ of
$W$, is a subalgebra and we may consider~$\exp(W_e)$.

In \cite{BZ} Y. A. Bahturin and M. V. Zaicev made the following
conjecture.

\begin{conjecture} \label{Conjecture BahZaicev}
$$
\exp(W) \leq |G|^2 \exp(W_e).
$$
\end{conjecture}

The conjecture was proved under certain natural conditions imposed
on the grading group $G$ and on the algebra~$W$ (see \cite{BZ}):

\begin{enumerate}

\item

The algebra~$W$ is finite-dimensional over a field of characteristic
zero and $W_{e}$ has polynomial growth.

\item

The algebra~$W$ is finitely generated semisimple over a field of
characteristic zero and $G$ is abelian.

\end{enumerate}

Our goal in this article is to prove the conjecture for affine
algebras where $G$ is an arbitrary finite group. For future
reference we record this in

\begin{theorem} \label{Main theorem} Let $G$ be a finite group. If $W$ is
an \textit{PI}-affine $G$-graded
algebra over a field $F$ of characteristic zero then $\exp(W)
\leq |G|^2 \exp(W_e).$
\end{theorem}

The first step (which is a key step) is to reduce the problem to
finite-dimensional $G$-graded algebras. This is obtained by invoking
the ``Representability Theorem" for $G$-graded affine algebras (see
\cite{AB}, Theorem 1.1). In order to state the theorem recall that
for a $G$-graded algebra one can define in a natural way $G$-graded
identities and can consider the $T$-ideal of $G$-graded identities
$\id_G(W)$ in the free $G$-graded algebra $F\langle X_G \rangle$.
The ungraded version of the theorem below was proved by Kemer (see
\cite{Kemer2}, \cite{Kemer3}).

\begin{theorem} [Representability of $G$-graded affine algebras]
Let $W$ be a PI-affine algebra over a field $F$ of characteristic zero and assume it is
$G$-graded. Then there exists a field extension~$K$ of $F$ and a
finite-dimensional algebra $A$ over $K$ such that
$\id_{G}(W)=\id_{G}(A)$. In particular $\id(W)=\id(A)$.
\end{theorem}

Thus, in order to prove the conjecture for \textit{PI}-affine
$G$-graded algebras it is sufficient to prove it for an arbitrary
finite-dimensional $G$-graded algebra $A$.

\begin{remark} We may assume that the field $F$ is algebraically closed.
Indeed, if $\overline{F}$ is any extension of $F$ (and in
particular its algebraic closure) and if
$\overline{A}=A\otimes_{F} \overline{F}$, then $\overline{A}$ is
$G$-graded by $\overline{A}_{g}=A_{g}\otimes \overline{F}$ and one
shows that the $G$-graded identities of $\overline{A}$ coincide
with the $G$-graded identities of $A$ (see \cite[Remark
1.5]{AB},\cite[Remark 1]{gz1}). It follows that the $n$-th
codimension $c_n(\overline{A})$ (resp. $c_n(\overline{A_e})$) over
$\overline{F}$ coincides with the $n$-th codimension $c_n(A)$
(resp. $c_n(A)$) over $F$.
\end{remark}

The problem can be reduced further, namely to the case where $G$ is
a simple group. We say that the conjecture holds for a group $G$ if
it holds for any algebra~$W$ which is $G$-graded.

\begin{lemma} \label{reduction}
Suppose that $H$ is a normal subgroup of $G$. If the conjecture
holds for $H$ and for $G/H$, then it holds for $G$. Hence, it is
sufficient to prove the conjecture for simple groups $G$.
\end{lemma}

\begin{proof}
Assume that $A$ is $G$-graded. Then it is $G/H$-graded in the
obvious way. Note that the $\overline{e}$-component of $A$
($\overline{e}\in G/H$ is the projection of the identity element of
$G$) coincides with $A_{H} = \oplus_{h\in H} A_{h}$ (with the
$G$-grading). Hence we have
$$
\exp(A) \leq |G/H|^{2}\exp(A_{H})
$$
and
$$
\exp(A_{H}) \leq |H|^{2} \exp(A_{e}).
$$
This gives
$$
\exp(A) \leq |G|^{2}\exp(A_{e})
$$
as desired.
\end{proof}

\begin{remark}
In fact, we will not use this reduction since our proof does not
become ``easier" if $G$ is assumed to be a simple group. However,
proving the conjecture for cyclic groups of prime order (and hence
for solvable groups) is substantially simpler.
\end{remark}

\end{section}

\begin{section}{Proofs}

In this section we prove Theorem \ref{Main theorem}. We assume (as
we may), that $A$ is a $G$-graded, finite-dimensional algebra over
$F$ and $F$ is algebraically closed.

We need to consider two decompositions of $A$. One, as an ordinary
algebra and the other as a $G$-graded algebra.

Let $J$ denote the Jacobson radical of $A$ and $\widetilde{A} \cong
A/J$ the semisimple quotient. By the Wedderburn-Malcev principal
theorem, there exists a semisimple subalgebra~$S$ of $A$, that is
isomorphic to $\widetilde{A}$. So we can write $A\cong S\oplus J$
where the isomorphism is (only) as vector spaces. Furthermore, the
subalgebra~$S$ may be decomposed into a direct product $S_1\times
\cdots \times S_d$ of simple algebras.

A similar decomposition holds if we consider $A$ as a $G$-graded
algebra. Indeed, it is well known that the Jacobson radical $J$ of
$A$ is $G$-graded (see \cite{CM}). Furthermore there exists a
semisimple subalgebra $S_G$ in $A$ that is
$G$-graded and is supplementary to $J$ as an $F$-vector space (see \cite{StVan}). The
algebra $S_G$ may be decomposed into the direct product of
$G$-simple algebras $(S_G)_1\times \cdots \times (S_G)_q$. Clearly,
we may decompose further each $G$-simple component into simple ones
and get a decomposition as above. We see that every simple component
$S_i$ is a constituent of precisely one $G$-simple component
$(S_G)_{j(i)}$. We refer to the $G$-simple algebra $(S_G)_{j(i)}$ as
the $G$-graded envelope (or just the envelope) of $S_i$.

Let us recall now how $\exp(A)$ is computed in terms of the
structure of $A$. Consider all possible \textit{nonzero} products in
$A$ of the form $S_{i_1}JS_{i_2}\cdots J S_{i_r}$. For such a
product consider the sum of the dimensions of the \textit{different}
simple algebras $S_i$ that are represented in the product. Then it
is known that $\exp(A)$ is equal to the maximal value of these
dimensions (see \cite{gz1}, \cite{gz2}). In particular this shows
that $\exp(A)$ is an integer.

Applying linearity we can replace each $S_i$ that appears in a
nonzero product $S_{i_1}JS_{i_2}\cdots J S_{i_r}$ by the
corresponding $G$-simple envelope and get a nonzero product
$(S_G)_{j(i_1)}J(S_G)_{j(i_2)}\cdots J (S_G){j(i_r)}$. Of course,
there may be repetitions among the $G$-simple components that appear
(even if there are no repetitions among the $S_i$'s) but clearly the
sum of the dimensions of the \textit{different} $G$-simple
components which appear in such nonzero product is at least the sum
of the dimensions of the~$S_i$'s. Let $\exp_{conj}^G(A)$ be the
maximal possible value obtained in this way (as sum of the
dimensions of the different $G$-simple components that appear in a
nonzero product).

\begin{remark}
We use the notation $\exp_{conj}^G(A)$ since it turns out to be equal to
$\exp^G(A)$, the $G$-graded exponent of the $G$-graded
algebra $A$, if $G$ is abelian and \textit{conjecturally}
if $G$ is arbitrary (see \cite{AGM}). In fact, for our purposes,
we will only need $\exp_{conj}^G(A)$ (and not $\exp^G(A)$). Nevertheless,
for completeness, we recall the definition
of $\exp^G(A)$. Let $\mathcal{W}_{G}=F\langle X,G \rangle$ be
the free $G$-graded algebra on a countable set, $\id_{G}(A)$ the
ideal of $G$-graded identities of $A$ and $F\langle X,G \rangle/
\id_{G}(A)$ the corresponding relatively free $G$-graded algebra.
For $n=1, 2,\ldots,$ we denote by $c^G_n(A)$ the dimension of the
subspace of $\mathcal{W}_{G}$ spanned by multilinear elements in $n$
($G$-graded) free generators. We call $\{c^G_n(A)\}_{n=1}^{\infty}$
the sequence of $G$-codimensions of $A$. The $G$-graded exponent of
$A$, is given by the limit $\lim_{n\to \infty}\sqrt[n]{c^G_n(A)}.$
\end{remark}

The following theorem is the main result of the paper.

\begin{theorem} \label{General Theorem}

With the above notation, $\exp_{conj}^G(A) \leq |G|^2 \exp(A_e)$.
\end{theorem}

Theorem \ref{Main theorem} follows from this since  $\exp(A) \leq
\exp_{conj}^G(A)$.

Our proof is based on a key result of Bahturin, Sehgal and Zaicev
(see \cite{BSZ}) in which they explicitly describe the structure of
$G$-graded, finite-dimensional $G$-simple algebras in terms of
\textit{fine} and \textit{elementary} gradings.

Before stating their result recall that a $G$-grading on a matrix
algebra $M_{r}(F)$ is said to be \textit{elementary} if there exists
an $r$-tuple $(g_{k_1},\dots,g_{k_r})\in G^{\times r}$ such that for
any $1\leq i,j \leq r$, the elementary matrix $e_{i,j}$ is
homogeneous of degree $g_{k_i}^{-1}g_{k_j}\in G$. A $G$-grading on an
algebra $A$ over $F$ is said to be \textit{fine} if each homogeneous
component~$A_{g}$, is of dimension $\leq 1$ (as an $F$-space).

\begin{theorem}[\cite{BSZ}]\label{BSZ}

Let $B$ be a $G$-simple algebra. Then there exists a subgroup $H$ of
$G$, a $2$-cocycle $f:H\times H\longrightarrow F^{*}$ where the
action of $H$ on $F$ is trivial, an integer $r$ and an $r$-tuple
$(g_{k_1},\ldots,g_{k_r})\in G^{\times r}$ such that $B$ is
$G$-graded isomorphic to $C=F^{f}H\otimes M_{r}(F)$, where $C_g=
\span_F\{u_h \otimes e_{i,j}: g=g_{k_i}^{-1}hg_{k_j}\}$. Here
$F^{f}H= \sum_{h\in H}Fu_{g}$ is the twisted group algebra of $H$
over $F$ with the $2$-cocycle $f$ and $e_{i,j}\in M_{r}(F)$ is the
$(i,j)$-elementary matrix.

In particular the idempotents $1 \otimes e_{i,i}$ as well as the
identity of $B$ are homogeneous of degree $e\in G$.
\end{theorem}

\begin{remark}
Note that the grading above induces $G$-gradings on the subalgebras
$F^{f}H\otimes F$ and $F\otimes M_{r}(F)$
which are \textit{fine} and \textit{elementary} respectively.
\end{remark}

For the clarity of the exposition it is convenient to ``place" the
$e$-component in the elementary grading of $M_{r}(F)$ as blocks
along the diagonal. Indeed reordering the $r$-tuple we can assume
the grading on $A\cong M_{r}$ is as follows:

Fix an ordering in the group $G$, $(g_1,\ldots,g_n)$. Fix also a
decomposition of $r$ into the sum of $n$ integers $r=r_1 + r_2 +
\ldots + r_n$. Note that some $r_i$ may be zero. Then we assign the
value $g_1=e$ to the first $r_1$ indices, the value $g_2$ to the
next $r_2$ indices and so on. Note that this determines a
$G$-grading on $M_r(F)$ where the $e$-blocks are concentrated along
the diagonal. We get $|G|$ $e$-blocks (again, we allow $0\times
0$-blocks) where the $i$-th block is of size $r_i$.

In the sequel, whenever we are given a $G$-simple algebra, we will
assume the grading is given as in Theorem \ref{BSZ} where the
elementary grading on $M_r(F)$ is as above.

Let us analyze the $e$-component of $B$ ($B$ is $G$-simple algebra
with the $G$-grading as in Theorem \ref{BSZ}). By the definition of
the $G$-grading, a basis element $u_{h}\otimes e_{i,j}$ is in the
$e$-component if and only if $e=g_{k_i}^{-1}hg_{k_j}$ or
equivalently $g_{k_i}=hg_{k_j}$. This means that $g_i$ and $g_j$
represent the same \textit{right} coset of $H$ in $G$. Consider the
family of right $H$-cosets in $G$. For convenience we order the
elements of $G$ (in the elementary grading of $M_r(F)$) according to
these equivalence classes. As a consequence we see that the
$e$-component of $B$ consists of an algebra which is isomorphic to
the direct product of matrix algebras $M_{t_1}(F) \times \dots
\times M_{t_l}(F)$ where $l=[G:H]$, $t_i$ is the number of elements
in the grading vector $(g_{k_1},\ldots,g_{k_r})\in G^{\times r}$
that belong to the $i$-th (right) $H$-coset, for $i=1,\ldots,l$
(note that $t_i$ may be zero). In particular every element of the
form $1\otimes e_{i,i}$ belongs to a unique simple component of
$B_{e}$. We refer to this component as the $e$-block of $B_{e}$
determined by the element $1 \otimes e_{i,i}$.

Let $A$ be a $G$-graded finite-dimensional algebra. Consider the
decomposition $A \cong S_{G} \oplus J$ where $S_{G}$ is a semisimple
algebra and $J$ is the Jacobson radical. Let $S_{G} \cong (S_G)_{1}
\times \cdots \times (S_G)_q$ be the decomposition of $S_{G}$ into
the direct product of its $G$-simple components. In view of the
above decomposition we will consider homogeneous elements that
belong to $J$ (radical elements) and basis elements of the form
$u_h\otimes e_{i,j} \in F^{f}H\otimes M_{r}(F)$ (semisimple
elements) where $F^{f}H\otimes M_{r}(F)$ is a $G$-simple component
of $S_G$.

In the following lemma we construct suitable monomials in $B$ which are
nonzero products of basis elements of the form $u_h\otimes e_{i,j}$.

\begin{lemma} [see proof of Lemma $5.2$ in \cite{AB}] \label{String}

Let $B=F^{f}H\otimes M_{r}(F)$ be a $G$-simple algebra ($G$-graded
as above).
\begin{enumerate}

\item

There exists a nonzero $G$-graded monomial that consists of the
product of all $r^2$ elementary matrices of $M_{r}(F)$ without
repetitions. Moreover, for any $k=1,\ldots,r$, there is such a
monomial which starts with the elementary matrix $e_{k,k}$ and
ends by $e_{i,k}$ for some $i\neq k$. Consequently we have a
nonzero monomial $Z$ which consists of all $r^2$ elements of the
form $1\otimes e_{i,j}$ in~$B$.

\item

For any $1 \leq k\leq r$ there exists a nonzero $G$-graded monomial
$\widehat{Z}$ which consists of basis elements of $F^{f}H\otimes
M_{r}(F)$ such that

\begin{enumerate}

\item

It starts with $1\otimes e_{k,k}$ and ends by an
element of the form $u_h\otimes e_{i,k}$

 \item

For any $1 \leq i\leq r$ the idempotent $1\otimes e_{i,i}$ appears
in $\widehat{Z}$.

\item

If $\widehat{Z}$ decomposes into the product $X\cdot1\otimes
e_{i,i}\cdot Y$ where $X$ has homogeneous degree $g\in G$, then
for any $\widehat{g}\in gg^{-1}_{k_{i}}Hg_{k_{i}}$ there exists a
decomposition $\widehat{Z}=X_{0}\cdot1\otimes e_{i,i}\cdot Y_{0}$
where $X_{0}$ has homogeneous degree~ $\widehat{g}$.
\item

The total value of $\widehat{Z}$ in $B$ is $1\otimes e_{k,k}$.

\end{enumerate}

\end{enumerate}

\end{lemma}

\begin{proof}
The first part is well known. In order to construct $\widehat{Z}$
we replace the idempotent $1\otimes e_{i,i}$ in $Z$ by the
monomial
$$u_{h_1}\otimes e_{i,i}\cdot z \cdot u_{h_2}\otimes e_{i,i}\cdot z \cdots z \cdot u_{h_r}\otimes e_{i,i}\cdot z$$
where $z=1\otimes e_{i,i}$, $h_{1}=1$ and

$$(h_{1}, h_{1}h_{2},\ldots, h_{1}h_{2}\cdots h_{r})$$
consist of all elements of the group $H$. This takes care of parts
(a) (b) and (c). The condition in part (d) is obtained by
multiplying (on the right) by an element of the form $\lambda
u_{h}\otimes e_{k,k}$.

\end{proof}

\begin{remark} \label{e-String}
Any semisimple element $u_{h}\otimes e_{i,j}$  may be multiplied
from left (resp. right) by $1\otimes e_{i,i}$ (resp. $1\otimes
e_{j,j}$) without changing its value and hence by the lemma, we can
insert a $G$-graded monomial $\widehat{Z}$ (as in Lemma
\ref{String}) from left (or right) and such that $\widehat{Z}\times
u_{h}\otimes e_{i,j} = u_{h}\otimes e_{i,j}$ (or $u_{h}\otimes
e_{i,j} \times \widehat{Z} = u_{h}\otimes e_{i,j}$). Furthermore
since the element $1\otimes e_{1,1}$ appears in the $G$-graded
monomial $\widehat{Z}$ we can insert an additional $G$-graded
monomial $E$ of the same form which starts with $1\otimes e_{1,1}$
and whose value is $1\otimes e_{1,1}$. Thus, we can decompose the
$G$-graded monomial $\widehat{Z}$ into a product $XY$ such that

$$
\widehat{Z}=XY=XEY
$$
where $E\in B_{e}$ is as above. Moreover for every semisimple
element $u_{h}\otimes e_{i,j}$ there exist $G$-graded monomials $X,
E, Y$ as above such that $XEY\times u_{h}\otimes e_{i,j} =
u_{h}\otimes e_{i,j}$ (similarly we could find $G$-graded monomials
which multiply $u_{h}\otimes e_{i,j}$ from the right).
\end{remark}

Let
$$\Lambda=z_{1}v_{1}z_{2}\cdots z_{n}v_{n}z_{n+1}$$
be a (nonzero) monomial in $A$ that realizes the value of
$\exp_{conj}^G(A)$. Here the $z$'s are semisimple elements and the
$v$'s are radical. Note that we may assume that the semisimple
elements belong to \textit{different} $G$-simple components. Indeed,
those that repeat may be ``swallowed'' by the radical elements.

Now, every semisimple element $z_{i}$ may be multiplied (say from
the right) by a $G$-graded monomial of semisimple elements (that
belong to the same $G$-simple component) of the form
$X_{i}E_{i}Y_{i}$ without changing the value of $\Lambda$.
Furthermore, the $G$-graded monomial $X_i$ may be ``swallowed'' by
the radical elements $v_{i}$ for $i=2,\ldots,n+1$ whereas $Y_j$ may
be ``swallowed'' by the radical elements $v_{j+1}$ for
$j=1,\ldots,n$. Finally, if we throw away $X_1$ and $Y_{n+1}$ from
the $G$-graded monomial we get a nonzero $G$-graded monomial

$$
\Omega = E_{1}v_{1}E_{2}\cdots E_{n}v_{n}E_{n+1}
$$
with the following properties:

\begin{enumerate}

\item

For $i=1,\ldots, n+1$, $E_{i}$ is a $G$-graded monomial whose
elements are in the $i$-th $G$-simple component (after renumbering).
In particular $n+1 \leq q$ where $q$ is the total number of
$G$-simple components of $A$.

\item

The (total) value of $E_{i}$ is $1\otimes e_{1,1}$ of the
corresponding $G$-simple component. Furthermore, the $G$-graded
monomial which corresponds to $E_{i}$ starts with $1\otimes
e_{1,1}$.

\item

The $G$-graded monomial $\Omega = E_{1}v_{1}E_{2}\cdots
E_{n}v_{n}E_{n+1}$ realizes the value of $\exp_{conj}^G(A)$. We may
refer to $\Omega$ as an element of $A$ and also as a $G$-graded
monomial whose elements are in $A$. We let $g_{0}\in G$ be the
homogeneous degree of $\Omega$.
\end{enumerate}

\begin{definition}
We say that a $G$-graded monomial $T=a_1\cdots a_n$ of homogeneous
elements ($a_{i}\in A$) has \textit{no proper} $g$-submonomial if either $n=1$ or else
there is no word of the form $T^{'}=a_1\cdots a_m$ for $m < n$ whose homogeneous degree is $g$.
\end{definition}

For any $g$ in $G$ we consider a possible decomposition of

$$
\Omega = X_{g}\Sigma_{1}\Sigma_{2}\cdots\Sigma_{d}Y_{g^{-1}g_{0}}
$$

where
\begin{enumerate}

\item

$X_{g}$ has homogeneous degree $g$ and has no proper $g$-submonomial.

\item

For $i=1,\ldots,d$, $\Sigma_{i}$ has homogeneous degree $e$ and has no proper $e$-submonomial.

\item

$Y_{g^{-1}g_{0}}$ has homogeneous degree $g^{-1}g_{0}$ and has no proper $e$-submonomial.
\end{enumerate}

\begin{remark}
\begin{enumerate}

\item
The decomposition above is maximal and unique in the sense that
there is no other decomposition with the same structure and the
same or more number of words (i.e. starting with a word of
homogeneous degree $g$, followed by words of homogeneous degree
$e$ and ending by a word of homogeneous degree $g^{-1}g_{0}$).
\item
Note that such a decomposition may not exist (indeed $X_{g}$ may not
exist). Note also that $d$ may be
zero, in which case we have $\Omega = X_{g}Y_{g^{-1}g_{0}}$. In what
follows there is no need to consider these cases separately
(including the case where the decomposition does not exist).
\end{enumerate}
\end{remark}

Given such a decomposition (for $g\in G$) we consider the semisimple
elements (idempotents) of the form $1\otimes e_{i,i}$ which can be
inserted in between any two adjacent words and giving nonzero
product. Note that these idempotents are uniquely determined. We
call these idempotents `` $e$-stops".

\begin{remark}

\begin{enumerate}
\item

Note that the different $e$-stops may or may not belong to the
same $G$-simple component.

\item

As mentioned above if an $e$-stop belongs to a $G$-simple component
$B$ then it determines uniquely an $e$-block of $B_{e}$.

\item
$e$-stops that belong to the same $G$-simple component determine the
same $e$-block of $B_{e}$.

\item

A $G$-simple component may not be represented by an $e$-stop.

\item
In case $X_{g}$ does not exists, the set of $e$-stops is empty.

\end{enumerate}
\end{remark}

Now for any $g\in G$, note that the $G$-graded monomial

$$
\Omega = X_{g}\Sigma_{1}\Sigma_{2}\cdots\Sigma_{d}Y_{g^{-1}g_{0}}
$$
gives rise to the \textit{nonzero} $G$-graded monomial

$$
(1\otimes e_{i_1,i_1})_{s_1} \Sigma_{1} (1\otimes e_{i_2,i_2})_{s_2}
\Sigma_{2} \cdots \Sigma_{d} (1\otimes
e_{i_{d+1},i_{d+1}})_{s_{d+1}}.
$$

Here the idempotent $(1\otimes e_{i_j,i_j})_{s_j}$,
$j=1,\ldots,d+1$, belongs to $(S_G)_{s_j}$, (the $s_j$-th $G$-simple
component in the decomposition of $A$). The reader should keep in
mind that all indices depend on $g\in G$.

Let us consider the above monomial as a monomial in $A_{e}$.
Clearly, since the monomial is nonzero, it provides a lower bound to
$\exp(A_e)$. Let us summarize the above considerations:

\begin{enumerate}
\item

The value of $\exp^{G}_{Conj}(A)$ is realized by the (nonzero)
$G$-graded monomial $$ \Omega = E_{1}v_{1}E_{2}\cdots
E_{n}v_{n}E_{n+1}.$$

\item

For any $g$ in $G$ we consider the decomposition

$$\Omega =
X_{g}\Sigma_{1}\Sigma_{2}\cdots\Sigma_{d}Y_{g^{-1}g_{0}}$$ and the
corresponding nonzero monomial in $A_{e}$

$$ (1\otimes e_{i_1,i_1})_{s_1} \Sigma_{1} (1\otimes
e_{i_2,i_2})_{s_2} \Sigma_{2} \cdots \Sigma_{d} (1\otimes
e_{i_{d+1},i_{d+1}})_{s_{d+1}}.$$

\end{enumerate}

\smallskip

As mentioned above each $e$-stop $(1\otimes e_{i_j,i_j})_{s_j}$
determines uniquely an $e$-block of the $G$-simple component
$(S_G)_{s_j}$. Note also (by the construction of the word $\Omega$)
that $e$-stops that belong to the same $G$-simple component are
adjacent to each other. In other words once we ``leave" a $G$-simple
component $(S_G)_{m}$ we do not return to it. Consequently, all
$e$-stops that belong to the same $G$-simple component determine the
same $e$-block.

As a result of this we have that each configurations (arising from
different elements in $G$) provides a lower bound to $\exp(A_{e})$.
More precisely:

\begin{corollary}

Denote by $(B_{g,1},\ldots, B_{g,n+1})$ the $e$-blocks (in the
different $G$-simple components) determined by $g\in G$ and let
$b_{g,i}^{2}=\dim_{F}(B_{g,i})$. Then $$\exp(A_{e})~\geq
b_{g,1}^{2} + \cdots + b_{g,n+1}^{2}.$$

\end{corollary}

Suppose now we start with a different element $\widehat{g}$ of
$G$. Consider the words $X_{\widehat{g}}$,
$Y_{\widehat{g}^{-1}g_{0}}$ and $\Sigma$'s. These determine
$e$-stops and consequently $e$-blocks. In the next lemma we
establish the relation between elements of $G$ which determine the
same $e$-block in a $G$-simple component. Let us denote by $H_{m}$
the group $H$ which appears in the $G$-graded structure of the
$G$-simple component $(S_{G})_{m}$ (see Theorem \ref{BSZ}).

\begin{lemma}
The number of elements in $G$ which determine the same (nonzero) $e$-block in $(S_{G})_{m}$ is precisely $ord(H_{m})$.
\end{lemma}

\begin{proof}
First note that each (nonzero) $e$-block in $(S_{G})_{m}$ is
determined by some $g\in G$. Suppose now that $g$ and
$\widehat{g}$ determine $e$-stops $1\otimes e_{i,i}$ and $1\otimes
e_{j,j}$ in $(S_{G})_{m}$ respectively. This means that
$\widehat{g}=gg_{k_i}^{-1}hg_{k_j}$ for some $h\in H_{m}$. But the
idempotents $1\otimes e_{i,i}$ and $1\otimes e_{j,j}$ determine
the same $e$-block in $(S_{G})_{m}$ (and hence $g$ and
$\widehat{g}$ determine the same $e$-block in $(S_{G})_{m}$) if
and only if $H_{m}g_{k_i}=H_{m}g_{k_j}$. Replacing $g_{k_j}$ by an
element of the form $h^{'}g_{k_i}$, $h^{'}\in H_{m}$, we see that
$g$ and $\widehat{g}$ determine the same $e$-block in
$(S_{G})_{m}$ if and only if they represent the same left coset of
$g_{k_i}^{-1}H_{m}g_{k_i}$ in $G$. In order to complete the proof
of the lemma it remains to show that if $g$ determines an $e$-stop
in $(S_{G})_{m}$ then all elements in $gg_{k_i}^{-1}H_{m}g_{k_i}$
determine $e$-stops in $(S_{G})_{m}$ as well. But this follows
from Lemma \ref{String} and so the lemma is proved.

\end{proof}

The main point of the lemma is the following

\begin{corollary}
Consider the different decompositions of $\Omega$

$$
\Omega = X_{g}\Sigma_{1}\Sigma_{2}\cdots\Sigma_{d}Y_{g^{-1}g_{0}}
$$
where $g$ runs over all elements of $G$. Then each (nonzero)
$e$-block in any $G$-simple component $(S_{G})_{m}$ will be
represented by the corresponding $e$-stops precisely $ord(H_{m})$
times. In the calculation below we include the $e$-blocks of
dimension $0$.
\end{corollary}

Let us show now that $\exp_{conj}^G(A)$ is bounded (from above) by
$|G|^2 \exp(A_e)$.
Suppose not. Then

$$
\exp_{conj}^G(A) > |G|^{2} \exp(A_e)
$$
and hence, for every $g\in G$ we have

$$\exp_{conj}^G(A) > |G|^{2} (b_{g,1}^{2} + \cdots + b_{g,n+1}^{2}).$$
Summing over all elements of $G$, we obtain

$$
|G|\exp_{conj}^G(A) > \sum_{g\in G} |G|^{2} (b_{g,1}^{2} + \cdots +
b_{g,n+1}^{2}).
$$
Recall now that every $G$-simple algebra $(S_G)_i$ has $r_i=[G:H_i]$
$e$-blocks (including the $0 \times 0$-blocks). In order to simplify
the notation we rename the $e$-blocks as follows: the $e$-blocks of
the $i$-th $G$-simple component $(S_G)_i$ will be denoted by
$B_{i,1},\ldots, B_{i,r_i}$ and $b_{i,j}^{2}=\dim_{F}(B_{i,j})$.

Now the left hand side yields

$$
|G|\exp_{conj}^G(A) = |G|\times \sum_{i=1}^{n+1}|H_i|(b_{i,1} + \cdots +
b_{i,r_i})^{2}
$$
whereas the right hand side yields
$$
|G|^{2} \sum_{i=1}^{n+1}|H_i|(b_{i,1}^{2}+ \cdots + b_{i, r_i}^{2}).
$$
In order to complete the proof we need to show that the inequality
$$
|G|\times \sum_{i=1}^{n+1}|H_i|(b_{i,1} + \cdots + b_{i,r_i})^{2}
> |G|^{2} \sum_{i=1}^{n+1}|H_i|(b_{i,1}^{2}+ \cdots + b_{i,
r_i}^{2})
$$
is not possible. Clearly, this will follow if for every $i$ we
have
$$
|G|\times |H_i|(b_{i,1} + \cdots + b_{i,r_i})^{2} \leq |G|^{2}
|H_i|(b_{i,1}^{2}+ \cdots + b_{i,r_i}^{2}).
$$

We simplify the notation once again by putting $H=H_i$, $r=r_i$ and
$b_{j}=b_{i,j}$ for $j=1,\ldots, r$. We therefore need to show that the
inequality
$$
(b_1 + \cdots + b_r)^{2} \leq |G|(b_{1}^2 + \cdots + b_{r}^2)
$$
holds. Clearly, it suffices to show that
$$
(b_1 + \cdots + b_r)^{2} \leq [G:H](b_{1}^2 + \cdots + b_{r}^2)=r(b_{1}^2 + \cdots + b_{r}^2)
$$
and this is clear since it is equivalent to
$$
\sum_{i<j}(b_i - b_j)^{2} \geq 0.
$$
This completes the proof of Theorem \ref{General Theorem} and hence of Theorem
\ref{Main theorem}.

\title{\textbf{Acknowledgment}: I would like to thank the referee for his comments on an earlier version
of this paper.}
\end{section}

\end{document}